\documentclass[12pt]{amsart} 

\oddsidemargin \evensidemargin
\newtheorem{theorem}{Theorem}[section]
\newtheorem{proposition}{Proposition}[section]

\newtheorem{lemma}[theorem]{Lemma}
\newtheorem{cor}[theorem]{Corollary}
\theoremstyle{definition}
\newtheorem{definition}[theorem]{Definition}

\newtheorem{procedure}[theorem]{Procedure}

\theoremstyle{remark}

\numberwithin{equation}{section}
\numberwithin{figure}{section}
\usepackage{color}
\usepackage{amsmath, latexsym, amssymb}
\usepackage{graphicx}
\usepackage{amsfonts}
\usepackage{subfigure}
\usepackage{amsthm}
\usepackage{amsmath}
\usepackage{graphics}
\usepackage{amssymb}
\usepackage{latexsym}
\usepackage{amscd}
\usepackage{color}
\usepackage{graphicx}
\usepackage{pict2e}

\newlength\cellsize 
\savebox2{%
\begin{picture}(15,15)
\put(0,0){\line(1,0){15}}
\put(0,0){\line(0,1){15}}
\put(15,0){\line(0,1){15}}
\put(0,15){\line(1,0){15}}
\end{picture}}
\newcommand\cellify[1]{\def\thearg{#1}\def\nothing{}%
\ifx\thearg\nothing
\vrule width0pt height\cellsize depth0pt\else
\hbox to 0pt{\usebox2\hss}\fi%
\vbox to 15\unitlength{
\vss
\hbox to 15\unitlength{\hss$#1$\hss}
\vss}}
\newcommand\tableau[1]{\vtop{\let\\=\cr
\setlength\baselineskip{-16000pt}
\setlength\lineskiplimit{16000pt}
\setlength\lineskip{0pt}
\halign{&\cellify{##}\cr#1\crcr}}}
\savebox3{%
\begin{picture}(15,15)
\put(0,0){\line(1,0){15}}
\put(0,0){\line(0,1){15}}
\put(15,0){\line(0,1){15}}
\put(0,15){\line(1,0){15}}
\end{picture}}
\newcommand\expath[1]{%
\hbox to 0pt{\usebox3\hss}%
\vbox to 15\unitlength{
\vss
\hbox to 15\unitlength{\hss$#1$\hss}
\vss}}

\DeclareMathOperator{\diagram}{dg}
\newcommand{\skydg}{\diagram'}
\newcommand{\augdg}{\widehat{\diagram}}

\DeclareMathOperator{\Des}{Des}
\DeclareMathOperator{\Inv}{Inv}

\DeclareMathOperator{\maj}{maj}

\DeclareMathOperator{\inv}{inv}
\DeclareMathOperator{\coinv}{coinv}

\DeclareMathOperator{\insr}{(\rho(F) \leftarrow k)}
\DeclareMathOperator{\insf}{(k \rightarrow F)}
\newcommand{\ns}{{\widehat E}_{\gamma}(x;0,0)}

\author[S. Mason]{S. Mason}
\thanks{Work partially supported by NSF postdoctoral research fellowship 
               DMS-0603351}
\address{Davidson College}

\email{samason@davidson.edu}
\urladdr{http://www.davidson.edu/math/mason}

\title[RSK Analogue]{A decomposition of Schur functions and an analogue of the Robinson-Schensted-Knuth Algorithm} 

\subjclass[2000]{Primary 05E05; Secondary 05E10}

\keywords{algebraic combinatorics, symmetric functions, Macdonald polynomials}

\begin{document}

\begin{abstract}

We exhibit a weight-preserving bijection between semi-standard Young tableaux and semi-skyline augmented fillings to provide a combinatorial proof that the Schur functions decompose into nonsymmetric functions indexed by compositions.  The insertion procedure involved in the proof leads to an analogue of the Robinson-Schensted-Knuth Algorithm for semi-skyline augmented fillings.  This procedure commutes with the Robinson-Schensted-Knuth Algorithm, and therefore retains many of its properties. \\

\end{abstract}

\maketitle

\section{Introduction}

Given a partition $\lambda=(\lambda_1,\lambda_2,\hdots)$, the {\it Schur function} $s_{\lambda}=s_{\lambda}(x)$ in the variables $x=(x_1,x_2,\hdots)$ is the formal power series
$s_{\lambda}=\sum_Tx^T$, where the sum is over all semi-standard Young tableaux $T$
of shape $\lambda$.  A {\it semi-standard Young tableau} is a diagram consisting of rows of squares such that the
$i^{th}$ row contains $\lambda_i$ squares, called {\it cells}.  This
diagram, called the Young (or Ferrers) diagram, is drawn in the first
quadrant, French style, as in~\cite{HHL05b}.  Each of these
cells is then assigned a positive integer entry so that the
entries are weakly increasing along rows and strictly
increasing along columns, where the rows are read from left to right and the columns are read from bottom to top.  The values assigned to the cells of $\lambda$
collectively form the multiset \{$1^{a_1},2^{a_2},\hdots,n^{a_n}$\}, for
some $n$, where $a_i$ is the number of times $i$ appears in $T$.  The weight of the semi-standard Young tableau (SSYT) $T$ is given by $x^T=\prod_{i=1}^{n}x_i^{a_i}$.  The sum of the weights of all SSYT of shape $\lambda$ is the Schur function $s_{\lambda}$.  See~\cite{Stanley2} for a more
detailed discussion of symmetric functions and Schur functions in
particular.

The Macdonald polynomials $\tilde{H}_{\lambda}(x;q,t)$ are a class
of functions symmetric in $x=(x_1,x_2,\hdots)$ with coefficients in $\mathbb{Q}(q,t)$.
Macdonald introduced them in \cite{Mac88} and conjectured that their
expansion in terms of Schur polynomials should have positive
coefficients.  A combinatorial formula for the Macdonald polynomials is proved by Haglund, Haiman, and Loehr in \cite{HHL05a}.  A corresponding combinatorial formula for the integral form Macdonald polynomials, $J_{\lambda}$, is obtained from the plethystic substitution which transforms the  $\tilde{H}_{\lambda}$ polynomials into integral form.

Building on this work, Haglund, Haiman, and Loehr \cite{HHL05c} derive a combinatorial formula for nonsymmetric Macdonald polynomials, which gives a new decomposition of $J_{\lambda}$ into nonsymmetric components.  Letting $q=t=\infty$ in this decomposition and using the fact that $J_{\lambda}(x;q,q)=s_{\lambda}(x)$ produces the identity 
\begin{equation}{\label{decomp}}
s_{\lambda}(x) = \sum_{\gamma \sim \lambda}{\widehat E}_{\gamma}(x;0,0),
\end{equation}
summed over all rearrangements $\gamma$ of the partition $\lambda$, denoted $\gamma \sim \lambda$.  (A {\it
rearrangement} of a partition $\lambda$ is an ordered sequence of non-negative integers (a weak composition), $\gamma$,
whose parts, when rearranged into weakly decreasing order, form the partition $\lambda$.)  Here ${\widehat E}_{\gamma}(x;0,0)$ is the specialization of the nonsymmetric Macdonald polynomial studied by Marshall \cite{M99}.  (We use the form ${\widehat E}$ instead of $E$ and replace $q$ and $t$ with $q^{-1}$ and $t^{-1}$ because the associated combinatorics is more elegant in this setting.)  The polynomial ${\widehat E}_{\gamma}(x;0,0)$ can be described combinatorially through objects called semi-skyline augmented fillings constructed using the statistics from \cite{HHL05c}. 

We provide a combinatorial proof of this decomposition of the Schur functions in Section \nolinebreak {\ref{comb}} by constructing a weight-preserving bijection $\Psi$ between semi-standard Young tableaux and semi-skyline augmented fillings.  We build a semi-skyline augmented filling through an insertion procedure similar to Schensted insertion.  This procedure is the fundamental operation in an analogue of the Robinson-Schensted-Knuth algorithm.

\begin{theorem}
{\label{rsk}}  There exists a bijection between $\mathbb{N}-$matrices with finite support and pairs $(F,G)$ of semi-skyline augmented fillings of compositions which rearrange the same partition.
\end{theorem}

This bijection is described in Section {\ref{rsk_sec}}.  We prove that it commutes with the Robinson-Schensted-Knuth (RSK) algorithm and retains the symmetry of the RSK algorithm plus many of its other properties.  We also describe the standardization of a semi-skyline augmented filling (SSAF) and prove that the RSK analogue commutes with standardization.  

\section{Combinatorial description of the functions $\ns$}

The functions $\ns$ are obtained from the integral form nonsymmetric Macdonald polynomials $\widehat{E_{\gamma}}(x;q,t)$ by letting $q$ and $t$ approach zero.  Haglund, Haiman, and Loehr provide a combinatorial formula for nonsymmetric Macdonald polynomials \cite{HHL05c} which can be specialized to obtain a combinatorial formula for nonsymmetric Schur functions.  Several definitions are necessary in order to describe this formula.

\subsection{Statistics on fillings}

Let $\gamma=(\gamma_1,\gamma_2,\hdots)$ be a weak composition of $n$ into $m$ parts.    The {\it column diagram of $\gamma$} is the figure $dg'(\gamma)$ consisting of $n$ {\it
cells} arranged into columns, as in \cite{HHL05c}.   The $i^{th}$ column contains $\gamma_i$ cells, and the number of cells in a column is called the {\it
height} of that column.   A cell $a$ in a column diagram is denoted $a=(i,j)$, where $i$ is the column index and $j$ is the row index.  

The {\it augmented diagram} of $\gamma$, defined by $\widehat{dg}(\gamma)=dg'(\gamma) \cup \{ (i,0) : 1 \le i \le m \}$ is the column diagram with $m$ extra cells adjoined in row $0$.  In this paper the adjoined row, called the {\it basement}, always contains the numbers $1$ through $m$ in strictly increasing order.  (We italicize the basement entries to avoid confusion.)

The augmented diagram for $\gamma=(0,2,0,3,1,2,0,0,1)$ is depicted below.
$$\augdg(\gamma)=\tableau{ & & & \mbox{} \\ & \mbox{} & & \mbox{} & & \mbox{} \\ & \mbox{} & & \mbox{} & \mbox{} & \mbox{} & & & \mbox{} \\ {\it 1} & {\it 2} & {\it 3} & {\it 4} & {\it 5} & {\it 6} & {\it 7} & {\it 8} & {\it 9}}$$

An {\it augmented filling,} $\widehat{\sigma}$, of an augmented diagram $\widehat{dg}(\gamma)$
is a function $\widehat{\sigma}: \widehat{dg}(\gamma) \rightarrow \mathbb{Z}_+$, which can be
pictured as an assignment of positive integer entries to the cells of
$\gamma$.  Let $\widehat{\sigma}(i)$ denote
the entry in the $i^{th}$ cell of the augmented diagram
encountered when $\widehat{dg}(\gamma)$ is read across rows from left to right, beginning at
the highest row and working down to the bottom row.  This ordering of the cells is called the {\it reading order}.  A cell $a=(i,j)$ precedes a cell $b=(i',j')$ in the reading order if either $j'<j$ or $j'=j$ and $i'>i$.  The reading word $\mathit{read}(\widehat{\sigma})$ is obtained by recording the non-basement entries in reading order.  The {\it content} of a filling $\widehat{\sigma}$ is the multiset of non-basement entries which appear in the filling.  (See Figure \ref{read example}.)

\begin{figure}[b]\label{read example}
$$\sigma=\tableau{ & & & \mbox{3} \\ & \mbox{2} & & \mbox{3} & & \mbox{1} \\ & \mbox{2} & & \mbox{4} & \mbox{5} & \mbox{6} & & & \mbox{9} \\ {\it 1} & {\it 2} & {\it 3} & {\it 4} & {\it 5} & {\it 6} & {\it 7} & {\it 8} & {\it 9}}$$
\caption{{\it read}$(\sigma)= 3 \; 2 \; 3 \; 1 \; 2 \; 4 \; 5 \; 6 \; 9$, {\it content}$(\sigma)=\{1,2^2, 3^2,4,5,6,9 \}$}
\end{figure}

A pair of cells $a$ and $b$ are called {\it attacking} if either $a$ and $b$ are in the same row or $a$ and $b$ are in adjacent rows, with the entry in the higher row strictly to the right of the entry in the lower row.  (That is, $a=(i_1,j_1)$ and $b=(i_2,j_2)$ are attacking if either $j_1=j_2$, or $j_2-j_1=1$ and $i_1 <i_2$, or $j_1-j_2=1$ and $i_2<i_1$.)   A {\it non-attacking filling} is a filling such that $\widehat{\sigma}(a) \not= \widehat{\sigma}(b)$ for every pair of attacking cells $a$ and $b$. 

Haglund, Haiman and Loehr introduce the statistics
$\Des(\widehat{\sigma})$ and $\Inv(\widehat{\sigma})$ to describe the nonsymmetric Macdonald polynomials.  As in \cite{HHL05b}, a {\it
descent} of $\widehat{\sigma}$ is a pair of entries $\widehat{\sigma}(a)>\widehat{\sigma}(b)$,
where the cell $a$ is directly above $b$.  (If $b=(i,j)$, then $a=(i,j+1)$.)  Call $\Des(\widehat{\sigma})=\{a \in
dg'(\gamma):\widehat{\sigma}(a)>\widehat{\sigma}(b)$ is a descent\} the {\it descent set} and define the {\it leg} of $u$ (denoted $l(u)$) to be the number of cells above $u$ in the column of $\widehat{dg}(\gamma)$ containing $u$.  Let $$\maj (\widehat{\sigma} ) = \sum _{u \in \Des (\widehat{\sigma} )} (l(u)+1).$$  Consider the cell $u=(i,j)$ and let $C_u$ be the column containing $u$.  The {\it arm} of $u$, denoted $a(u)$, is the number of cells to the right of $u$ in row $i$ appearing in columns whose height is weakly less than the height, $h$, of $C_u$ plus the number of cells to the left of $u$ in row $i-1$ appearing in columns whose height is strictly less than $h$.

Let $a_1=(i_1,j_1), a_2 = (i_2,j_2),$ and $a_3 =(i_3,j_3)$ be three cells in $\widehat{dg}(\gamma)$ such that column $i_1$ is taller than or equal in height to column $i_2$.  If $j_1=j_2$, $j_1-j_3=1$, and $i_1=i_3$, then $a_1,a_2,$ and $a_3$ are said to form a {\it type $A$ triple}, as depicted below.
$$\tableau{a_1 \\ a_3} \cdots \tableau{a_2}$$

Define for $x,y \in \mathbb{Z}_+$
\[I(x,y) = \left\{ \begin{array}{ll}
1  & \mbox{if $x>y$} \\
0  & \mbox{if $x \le y$}
\end{array}
\right. . \] 

Let $\widehat{\sigma}$ be an augmented filling and let $\{ \widehat{\sigma}(a_1),\widehat{\sigma}(a_2),\widehat{\sigma}(a_3) \}$ be the entries of
$\widehat{\sigma}$ in the cells $\{a_1,a_2,a_3\}$, respectively, of a type $A$ triple.  The triple $\{a_1,a_2,a_3\}$ is called a {\it type $A$ inversion triple} if and only
if $I(\widehat{\sigma}(a_1),\widehat{\sigma}(a_2))+I(\widehat{\sigma}(a_2),\widehat{\sigma}(a_3))-I(\widehat{\sigma}(a_1),\widehat{\sigma}(a_3))=1$.

Similarly, consider three cells $\{a_1=(i_1,j_1),a_2=(i_2,j_2),a_3=(i_3,j_3) \} \in \lambda$ such that column $i_2$ is strictly taller then column $i_1$.  The cells $\{a_1,a_2,a_3 \}$ are said to form a {\it type $B$ triple} if $j_1=j_2$, $i_2=i_3$, and $j_3-j_2=1$, as shown below.

$$\tableau{ \\ a_1} \hdots  \tableau{a_3 \\ a_2}$$

Let $\widehat{\sigma}$ be an augmented filling and let $\{ \widehat{\sigma}(a_1), \widehat{\sigma}(a_2), \widehat{\sigma}(a_3) \}$ be the entries of
$\widehat{\sigma}$ in the cells $\{a_1,a_2,a_3 \}$ of a type $B$ triple.  The triple $\{ a_1,a_2,a_3 \}$ is called a {\it type $B$ inversion triple} if and
only if $I(\widehat{\sigma}(a_3),\widehat{\sigma}(a_1))+I(\widehat{\sigma}(a_1),\widehat{\sigma}(a_2))-I(\widehat{\sigma}(a_3),\widehat{\sigma}(a_2))=1$.  

Let $\inv (\widehat{\sigma })$ be the number of type $A$ inversion triples plus the number of type $B$ inversion triples.  Let $\coinv (\widehat{\sigma })$ be the number of type $A$ or $B$ triples which are not inversion triples.  These two statistics appear in the formula for nonsymmetric Macdonald polynomials $E_{\gamma}$ given by Haglund, Haiman, and Loehr \cite{HHL05c}.  For our purposes it will be more convenient to work with the integral form nonsymmetric Macdonald polynomials ${\widehat E}_{\gamma}$ studied by Marshall \cite{M99}, which are related to the $E_{\gamma}$ via the equation
${\widehat E}_{\gamma_p, \hdots ,\gamma_1}(x_n,\ldots ,x_1;1/q,1/t)=E_{\gamma}(x_1, x_2, \ldots, x_n ; q,t)$.

\begin{theorem}{\cite{HHL05c}}
The nonsymmetric Macdonald polynomials ${\widehat E}_{\gamma}$ are given by the
formula
\begin{equation}
{\widehat E}_{\gamma }(x;q,t) = \sum _{\substack{\sigma \colon \gamma \rightarrow [n]\\
\text{non-attacking}}} x^{\sigma } q^{\maj (\widehat{\sigma })} t^{\coinv
(\widehat{\sigma })} \prod _{\substack{u\in \skydg (\gamma )\\
\widehat{\sigma }(u)\not =\widehat{\sigma }(d(u))}}
\frac{1-t}{1-q^{l(u)+1}t^{a(u)+1}},
\end{equation}
where $x^{\sigma } = \prod _{u\in \skydg (\mu )} x_{\sigma (u)}$.
\end{theorem}

We are concerned only with the polynomial ${\widehat E}_{\gamma}(x; 0,0)$, so setting $q$ and $t$ equal to zero produces the polynomial
\begin{eqnarray*}
{\widehat E}_{\gamma }(x;0,0) & = & \sum _{\substack{\sigma \colon \gamma \rightarrow [n]\\
\text{non-attacking}}} x^{\sigma } 0^{\maj (\widehat{\sigma })} 0^{\coinv
(\widehat{\sigma })} \prod _{\substack{u\in \skydg (\gamma )\\
\widehat{\sigma }(u)\not =\widehat{\sigma }(d(u))}}
\frac{1-0}{1-0^{l(u)+1}0^{a(u)+1}}\\
& = & \sum_{\substack{\sigma \colon \gamma \rightarrow [n]\\  \text{non-attacking}\\ \maj (\widehat{\sigma })=\coinv(\widehat{\sigma })=0}} x^{\sigma}
\end{eqnarray*}

A non-attacking filling $F$ satisfying $\maj (F)=\coinv(F)=0$ is called a {\it semi-skyline augmented filling, SSAF}.  In the following section, we provide a simpler definition of a semi-skyline augmented filling. 

\subsection{Semi-skyline augmented fillings}

Let $F$ be a semi-skyline augmented filling.  The condition $\maj(F)=0$ implies that $$\maj (F) = \sum _{a \in \Des (F)} (l(a)+1)=0.$$  Since $\maj (F)$ is increased by at least $1$ for each element in the descent set of $F$, we must have $\Des(F)=\emptyset$.  The condition $\coinv(F)=0$ implies that every triple of cells in $F$ must be an inversion triple.  Therefore a semi-skyline augmented filling is a non-attacking filling with no descents such that every triple is an inversion triple.

The following two lemmas demonstrate that any filling satisfying these descent and inversion conditions must be a non-attacking filling.

\begin{lemma}
{\label{equal}}  Let $F$ be a descentless augmented filling such that every triple of $F$ is an inversion triple.  If the cells $a$ and $b$ are in the same row of $F$, then $F(a) \not= F(b)$.
\end{lemma}

\begin{proof}
Suppose $a$ and $b$ are in the same row of $F$.  We may assume that $a$ is to the left of $b$.  Assume first that the column containing $a$ is weakly taller than the column containing $b$.  The cell $a$ is directly on top of some cell $c$, so $\{ a, b, c \}$ is a type $A$ triple as depicted below.
$$\tableau{{\bf a} \\ c } \hdots \tableau{{\bf b}}$$
The triple must be a Type $A$ inversion triple, so $I(F(a),F(b))+I(F(b),F(c))-I(F(a),F(c))=1$.  We know that $F$ contains no descents, so $I(F(a),F(c))=0$.  Therefore $I(F(a),F(b))+I(F(b),F(c))=1$ implies that either $F(a) > F(b)$ and $F(b) \le F(c)$ or $F(a) \le F(b)$ and $F(b) > F(c)$.  If $F(a) \le F(b)$, then $F(b) > F(c)$ implies that $F(a) < F(b)$, for otherwise $F(a)=F(b) > F(c)$, which contradicts the fact that $F(a) \le F(c)$.

Next suppose that the column containing $b$ is strictly taller than the column containing $a$.  There must be a cell $d$ on top of $b$ such that $\{a,d,b\}$ is a type $B$ triple as depicted below.  
$$\tableau{ \\ {\bf a}} \hdots \tableau{ d \\  {\bf b}}$$
There are no descents in $F$, so $F(d) \le F(b)$.  Every triple of $F$ is an inversion triple, so $I(F(d),F(a))+I(F(a),F(b))-I(F(d),F(b))=1$ implies that $I(F(d),F(a))+I(F(a),F(b))-0=1$.  We only need to consider the situation in which $F(a) \le F(b)$, which implies that $I(F(a),F(b))=0$.  This means that $I(F(d),F(a))=1$, so $F(a) < F(d)$.  But $F(d) \le F(b)$, so $F(a) < F(d) \le F(b)$ and therefore $F(a) < F(b)$.  So $F(a) \not= F(b)$.
\end{proof}

\begin{lemma}
{\label{attack}}  
Let $F$ be a descentless augmented filling such that every triple of $F$ is an inversion triple. 
For each pair of cells $a$ and $b$ in $F$, with $a$ to the left of $b$ in the row immediately below $b$, we have $F(a) \not= F(b)$.
\end{lemma}

\begin{proof}
Consider two cells $a$ and $b$ in the augmented filling situated as described.  There exists a cell $d$ immediately below $b$ and possibly a cell $c$ immediately above $a$ as depicted below.
$$\tableau{c \\ {\bf a} } \hdots  \tableau{{\bf b} \\ d}$$
If the column containing $a$ is taller than or equal to the column containing $b$, then $a$ lies directly below the cell $c$ which must have $F(c) \le F(a)$.  Since the triple $\{c,b,a\}$ is a type $A$ triple, it must satisfy $I(F(c),F(b))+I(F(b),F(a))-I(F(c),F(a))=1$.  Since $F(c) \le F(a)$, we have $I(F(c),F(a))=0$ and hence $I(F(c),F(b))+I(F(b),F(a))=1$.  We only need to consider the situation in which $F(b) \le F(a)$, which implies that $I(F(b),F(a))=0$.  In this case, $I(F(c),F(b))=1$, so $F(b) < F(c) \le F(a)$.  This means that $F(a) > F(b)$ and hence $F(a) \not= F(b)$.

If the column containing $b$ is strictly taller than the column containing $a$, the triple $\{a,b,d\}$ is a type $B$ triple and must satisfy $I(F(b),F(a))+I(F(a),F(d))-I(F(b),F(d))=1$.  There are no descents in $F$, so $F(b) \le F(d)$, which means that $I(F(b),F(d))=0$.  Therefore either $I(F(b),F(a))=1$ or $I(F(a),F(d))=1$.  If $I(F(b),F(a))=1$, then $F(b) > F(a)$ and we are done.  Assume $I(F(b),F(a))=0$ and $I(F(a),F(d))=1$.  This means that $F(a) > F(d) \ge F(b)$.  So $F(a) > F(b)$ and hence $F(a) \not= F(b)$.
\end{proof}

\begin{cor}{\label{noattacks}}
The descent and inversion conditions used to describe the semi-skyline augmented fillings are enough to guarantee that the filling is non-attacking.
\end{cor}

Corollary {\ref{noattacks}} follows immediately from Lemmas \ref{equal} and \ref{attack}.  It allows us to reformulate the combinatorial interpretation of $\widehat{E}_{\gamma}(x;0,0)$ as follows.

\begin{definition}
Let $\gamma$ be a weak composition of $n$ into $m$ parts (where $m \in \mathbb{Z}^+ \cup \{ \infty \}$).  The polynomial
$\widehat{E}_{\gamma}(x;0,0)$ in the variables $x=(x_1,x_2,...,x_k)$
is the formal power series $$\widehat{E}_{\gamma}(x;0,0)=\sum_{F \in SSAF(\gamma)} x^F,$$  where $SSAF(\gamma)$ is the set of all descentless fillings of $\widehat{dg}(\gamma)$ in which every triple is an inversion triple. 
\end{definition}

\begin{figure}
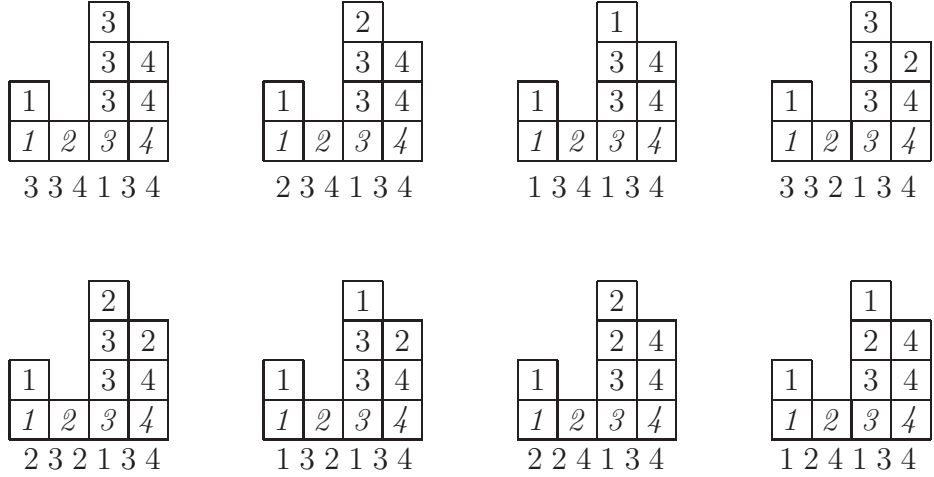

{\label{ns}}
$$\tableau{ & & 3 \\ & & 3 & 4 \\ 1 & & 3 & 4 \\ {\it 1} & {\it 2} & {\it 3} & {\it 4}} \hspace*{.5in} \tableau{ & & 2 \\ & & 3 & 4 \\ 1 & & 3 & 4 \\ {\it 1} & {\it 2} & {\it 3} & {\it 4}} \hspace*{.5in} \tableau{ & & 1 \\ & & 3 & 4 \\ 1 & & 3 & 4 \\ {\it 1} & {\it 2} & {\it 3} & {\it 4}} \hspace*{.5in} \tableau{ & & 3 \\ & & 3 & 2 \\ 1 & & 3 & 4 \\ {\it 1} & {\it 2} & {\it 3} & {\it 4}}$$
$$ 3 \; 3 \; 4 \; 1 \; 3 \; 4 \hspace*{.6in} 2 \; 3 \; 4 \; 1 \; 3 \; 4 \hspace*{.6in} 1 \; 3 \; 4 \; 1 \; 3 \; 4 \hspace*{.6in} 3 \; 3 \; 2 \; 1 \; 3 \; 4$$
\vspace*{.2in}
$$\tableau{ & & 2 \\ & & 3 & 2 \\ 1 & & 3 & 4 \\ {\it 1} & {\it 2} & {\it 3} & {\it 4}} \hspace*{.5in} \tableau{ & & 1 \\ & & 3 & 2 \\ 1 & & 3 & 4 \\ {\it 1} & {\it 2} & {\it 3} & {\it 4}} \hspace*{.5in} \tableau{ & & 2 \\ & & 2 & 4 \\ 1 & & 3 & 4 \\ {\it 1} & {\it 2} & {\it 3} & {\it 4}} \hspace*{.5in} \tableau{ & & 1 \\ & & 2 & 4 \\ 1 & & 3 & 4 \\ {\it 1} & {\it 2} & {\it 3} & {\it 4}}$$
$$ 2 \; 3 \; 2 \; 1 \; 3 \; 4 \hspace*{.6in} 1 \; 3 \; 2 \; 1 \; 3 \; 4 \hspace*{.6in} 2 \; 2 \; 4 \; 1 \; 3 \; 4 \hspace*{.6in} 1 \; 2 \; 4 \; 1 \; 3 \; 4$$
\caption{All possible SSAFs of shape $(1,0,3,2)$ with the reading word listed below each filling.  $\widehat{E}_{(1,0,3,2)}(x;0,0)=x_1x_3^3x_4^2+x_1x_2x_3^2x_4^2+ x_1^2x_3^2x_4^2+x_1x_2x_3^3x_4+x_1x_2^2x_3^2x_4 + x_1^2x_2x_3^2x_4 + x_1x_2^2x_3x_4^2 + x_1^2x_2x_3x_4^2$.}
\end{figure}

The combinatorial interpretation of $\widehat{E}_{(1,0,3,2)}(x; 0,0)$ is depicted in Figure {\ref{ns}}.  Notice that the first row of the $i^{th}$ column contains the entry $i$, since any smaller entry $j$ would attack the basement entry $j$ and any larger entry would create a descent.  The following lemma characterizes the type $B$ triples in an SSAF.

\begin{lemma}
{\label{typeb}}  If $\{a,b,c\}$ is a type $B$ triple in a semi-skyline augmented filling $F$ with $a$ and $c$ in the same row and $b$ directly above $c$, then $F(a) < F(c)$.
\end{lemma}

\begin{proof}
Let $\{a,b,c\}$ be a type $B$ triple situated as described:
$$\tableau{\\ a} \hdots \tableau{ b \\ c}$$

Lemma {\ref{equal}} implies that $F(a) \not= F(c)$.  If $F(a) > F(c)$, then consider the cells $e$ and $g$ immediately below $a$ and $c$ as shown.
$$\tableau{a \\e} \hdots \tableau{ c \\ g}$$
The descent condition implies that $F(a) \le F(e),$ which means that $I(F(c),F(e))=0$, since $F(c) < F(a) \le F(e)$.  Lemma \ref{attack} therefore implies that $F(c) < F(e)$. Since $\{e,c,g\}$ form a type $B$ inversion triple, we must have $I(F(c),F(e))+I(F(e),F(g))-I(F(c),F(g))=1$.  Therefore $0+I(F(e),F(g))-I(F(c),F(g))=1$ implies that $I(F(e),F(g))=1$.   By definition this means that $F(e) > F(g)$.  We now have the same situation as above, but one row lower in the diagram.  Repeating the argument eventually implies that the basement entry in the column containing $a$ is greater than the basement entry in the column containing $b$, a contradiction.  Therefore $F(a) < F(c)$.
\end{proof}

Given any type $B$ triple $\{a,b,c\}$ as described in Lemma {\ref{typeb}}, the descent condition implies that $F(b) \le F(c)$.  If $F(a) \ge F(b)$, then Lemma {\ref{typeb}} and the fact that the filling is non-attacking implies that $F(b) < F(a) < F(c)$.  This means that $I(F(b),F(a))+I(F(a),F(c))-I(F(b),F(c))=0+ 0-0=0$, so $\{ b, a, c \}$ is a type $B$ non-inversion triple.  Thus $F(a) < F(b) \le F(c)$.  Therefore Lemma {\ref{typeb}} completely characterizes the relative values of the cells in a type $B$ inversion triple.

\begin{lemma}{\label{typea}}
Consider two cells $a_1=(i_1,j_1)$ and $a_2=(i_2,j_2)$ in an augmented diagram $\augdg(\gamma)$ such that $i_2=i_1-1$ and $j_2 < j_1$.  Let $a_3=(i_3,j_3)$ be the cell directly on top of $a_2$.  (Then $j_3 = j_2+1$ and $i_3=i_2$.)  If $F$ is a semi-skyline augmented filling and $F(a_1) \le F(a_2)$, then $F(a_3) > F(a_1)$.
\end{lemma}

\begin{proof}
The cells $a_1,a_2,a_3$ form a type $A$ triple.  The filling $F$ must not contain any descents, so $F(a_3) \le F(a_2)$.  If $F(a_3) \le F(a_1)$, then $I(F(a_3),F(a_1))+I(F(a_1),F(a_2))-I(F(a_3),F(a_2))= 0+0-0=0 \not=1$.  However, if $F(a_3) > F(a_1)$, then $I(F(a_3),F(a_1))+I(F(a_1),F(a_2))-I(F(a_3),F(a_2))=1+0-0=1$, as required.  Therefore $F(a_3) > F(a_1)$.
\end{proof}

Lemmas {\ref{equal}}, {\ref{attack}}, {\ref{typeb}}, and {\ref{typea}}  provide several conditions on the cells in a semi-skyline augmented filling. These conditions are used to describe the bijection between semi-standard Young tableaux and semi-skyline augmented fillings and to prove Theorem {\ref{rsk}}.

\section{Two equivalent bijections between SSAFs and SSYTs}{\label{comb}}

Recall that we seek a weight-preserving bijection between semi-standard Young tableaux of shape $\lambda$ and semi-skyline augmented fillings whose shape rearranges $\lambda$.  Such a bijection provides the first combinatorial proof that
\begin{equation}
\sum_{\gamma \sim \lambda}\widehat{E}_{\gamma}(x;0,0)=s_{\lambda}(x),
\label{eq:nonsym}
\end{equation}
where the sum is over all infinite weak compositions $\gamma$ which rearrange $\lambda$, denoted $\gamma \sim \lambda$.  Set $x_i=0$ for $i>n$ to obtain the equation $$\sum_{\gamma \sim \lambda}\widehat{E}_{\gamma}(x_1,x_2,\hdots, x_n;0,0)=s_{\lambda}(x_1,x_2,\hdots, x_n).$$

\subsection{A weight-preserving bijection between SSAFs and reverse SSYTs}{\label{rev}}

Assume that $F$ is a semi-skyline augmented filling of a weak composition $\gamma$ of $n$ into infinitely many parts whose largest part is equal to $m$.  Let $R_i$ be the set of entries in the $i^{th}$ row of $F$.  The collection $\{ R_i \}_{i \in [m]}$ of all sets of row entries is called the {\it row set} of $F$.  

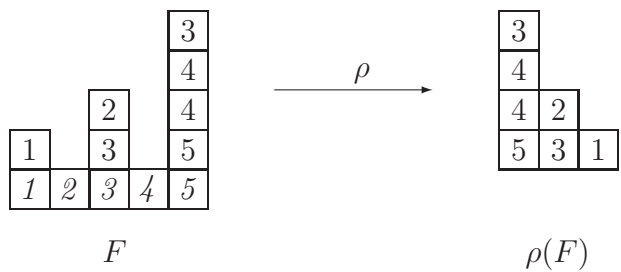
\begin{figure}{\label{rho}}
\begin{center}
\begin{picture}(230,85)
\put(35,0){$F$}
\put(0,80){$$\tableau{ & & & & 3 \\ & & & & 4 \\ & & 2 & & 4 \\ 1 & & 3 & & 5 \\ {\it 1} & {\it 2} & {\it 3} & {\it 4} & {\it 5}}$$}
\put(100,65){\vector(1,0){60}}
\put(130,70){$\rho$}
\put(185,80){$$\tableau{3 \\ 4 \\ 4 & 2 \\ 5 & 3 & 1}$$}
\put(195,0){$\rho(F)$}
\end{picture}
\end{center}
\caption{The map $\rho(F)$}
\end{figure}

Recall that a {\it reverse semi-standard Young tableau} is a filling of a partition shape with positive integer values whose entries are strictly decreasing along rows and weakly decreasing along columns, where the rows are read from left to right and the columns are read from bottom to top.  We construct a reverse semi-standard Young tableau $\rho(F)$ by placing the entries from the set $R_i$ into the $i^{th}$ row of $\rho(F)$ in decreasing order from left to right.  (See Figure \ref{rho}.)  The entries in $R_i$ are distinct since $F$ is a non-attacking filling.  Therefore the row entries of $\rho(F)$ are strictly decreasing.  

To see that the column entries of $\rho(F)$ are weakly decreasing, notice that any element in $R_i$ lies immediately above an element of $R_{i-1}$ in $F$.  No descents occur in $F$, so the $k^{th}$ largest element, $\alpha$, in $R_i$ must be less than or equal to at least $k$ of the elements in $R_{i-1}$.  The entry $\alpha$ appears in $\rho(F)$ immediately on top of the $k^{th}$ largest element, $\beta$, of $R_{i-1}$.  At least $k$ elements in $R_{i-1}$ are greater than or equal to $\alpha$, so $\beta \ge \alpha$.  Therefore the column entries of $\rho(F)$ are weakly decreasing and $\rho(F)$ is indeed a reverse semi-standard Young tableau.

The inverse of the map $\rho$ sends a reverse semi-standard Young tableau $P$ to a semi-skyline augmented filling $\rho^{-1}(P)$ as follows.  Assume that the first $i$ rows of $P$, denoted $\{P_1,P_2, \hdots P_i\}$ have been mapped to a semi-skyline augmented filling.  Consider the largest element, $\alpha_1$, in the $(i+1)^{th}$ row $P_{i+1}$.   There exists an element greater than or equal to $\alpha_1$ in the $i^{th}$ row of the SSAF since $\alpha_1$ is immediately above such an element in $P$.  Place $\alpha_1$ on top of the leftmost such element.  

Assume that the largest $k-1$ entries in $P_{i+1}$ have been placed into the SSAF.  The $k^{th}$ largest element, $\alpha_k$, of $R_{i+1}$ is then placed into the SSAF.  There are at least $k$ elements of $P_i$ which are greater than or equal to $\alpha_k$, so at most $k-1$ of these are already directly beneath an element of $P_{i+1}$.  Place $\alpha_k$ on top of the leftmost entry $\beta$ in row $k-1$ such that $\beta \ge \alpha_k$ and the cell immediately above $\beta$ is empty.  Again, such a $\beta$ exists since the number of entries in row $P_{i+1}$ which are greater than or equal to $\alpha_k$ is less than or equal to the number of entries in row $P_i$ which are greater than or equal to $\alpha_k$.  Continue this procedure until all entries in $P_{i+1}$ have been mapped into the $(i+1)^{th}$ row and then repeat for the remaining rows of $P$ to obtain the semi-skyline augmented filling $\rho^{-1}(P)$.  (See Figure 3.2.)

We must prove that $\rho^{-1}(P)$ is a semi-skyline augmented filling, and that this is the only SSAF with row entries $\{P_i\}$.

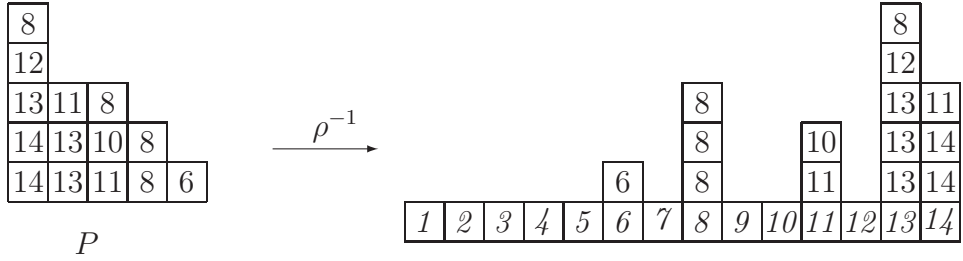
\begin{figure}{\label{rhoinv}}
\begin{center}
\begin{picture}(350,105)
\put(0,80){$$\tableau{8 \\ 12 \\ 13 & 11 & 8 \\ 14 & 13 & 10 & 8 \\ 14 & 13 & 11 & 8 & 6}$$}
\put(25,0){$P$}
\put(100,40){\vector(1,0){40}}
\put(115,45){$\rho^{-1}$}
\put(150,80){$$\tableau{ & & & & & & & & & & & & 8 \\ & & & & & & & & & & & & 12 \\ & & & & & & & 8 & & & & & 13 & 11 \\ & & & & & & & 8 & & & 10 & & 13 & 14 \\ & & & & & 6 & & 8 & & & 11 & & 13 & 14 \\ {\it 1} & {\it 2} & {\it 3} & {\it 4} & {\it 5} & {\it 6} & {\it 7} & {\it 8} & {\it 9} & {\it 10} & {\it 11} & {\it 12} & {\it 13} & {\it 14} }$$}
\end{picture}
\end{center}
\caption{The map $\rho^{-1}(P)$}
\end{figure}

\begin{lemma}{\label{inv}}
The filling $\rho^{-1}(P)=F'$ is a semi-skyline augmented filling.
\end{lemma}

\begin{proof}
By construction, the filling has no descents. Therefore, we must show that all triples are inversion triples.  First consider a type $A$ triple of cells, $\{a,b,c\}$ where $a$ and $b$ are cells in the $i^{th}$ row and $c$ is the cell immediately below $a$ as depicted below.  $$\tableau{a \\ c} \hdots \tableau{b}$$  The column containing $a$ and $c$ lies to the left of the column containing $b$ and therefore must be weakly taller than the column containing $b$.

The filling is constructed in such a way that $F'(a) \le F'(c)$.  Therefore $I(F'(a),F'(c))=0$.  We want to show $1=I(F'(a),F'(b))+I(F'(b),F'(c))-I(F'(a),F'(c))=I(F'(a),F'(b))+I(F'(b),F'(c))-0.$  If $I(F'(a),F'(b))=0,$ then $F'(a) < F'(b)$, since $a$ and $b$ are in the same row of $F'$. Since $F'(b) > F'(a)$, the entry $F'(b)$ was inserted into $F'$ before the entry $F'(a)$.  Since $F'(b)$ was not placed on top of $c$, the entry $F'(c)$ must be less than $F'(b)$.  Therefore $I(F'(b),F'(c))=1$.  So $I(F'(a),F'(b))+I(F'(b),F'(c))-I(F'(a),F'(c))=0+1-0=1$ and $\{a,b,c\}$ is a type $A$ inversion triple.

If $I(F'(a),F'(b))=1,$ then $F'(a)>F'(b)$.  This implies that $F'(b) < F'(c)$, since $F'(a) \le F'(c)$.  So $I(F'(a),F'(b))+I(F'(b),F'(c))-I(F'(a),F'(c))=1+0-0=1$.  So every type $A$ triple in $F'$ is a type $A$ inversion triple.

Next consider a type $B$ triple $\{a,b,c\}$ such that $a$ and $b$ are in the $i^{th}$ row and $a$ is to the left of $b$ as shown below.  $$\tableau{ \\ a} \hdots \tableau{c \\ b}$$  Let $c$ be the cell immediately on top of $b$. Then the column containing $b$ and $c$ is strictly taller than the column containing $a$.  We have the equation $I(F'(c),F'(a))+I(F'(a),F'(b))-I(F'(c),F'(b))=I(F'(c),F'(a))+I(F'(a),F'(b))-0$ since $F'(c) \le F'(b)$.

If $I(F'(a),F'(b))=1$, then $F'(a)>F'(b)$.  Since $F'(c)$ was not placed on top of $F'(a)$, there must be an entry $F'(d)$ greater than $F'(c)$ which was already on top of $F'(a)$ when $F'(c)$ was inserted.  If this is the case, then the same situation occurs in row $i+1$.  Repeat this argument until a row $r$ is reached such that the column containing $b$ and $c$ contains an entry $\alpha$ in row $r$ but the column containing $a$ does not.  (Such a row must exist since the triple is a type $B$ triple.)  Then the entry $\alpha$ would have been placed on the column containing $a$, which is a contradiction.  Therefore $I(F'(a),F'(b))=0$.

Since $I(F'(a),F'(b))=0$ implies $F'(a) \le F'(b)$ and the row entries of $P$ are distinct, we have $F'(a) < F'(b)$. Since $F'(c)$ was not placed on top of $F'(a)$, either $F'(c) > F'(a)$ or there was already an entry larger than $F'(c)$ on top of $a$.  If there was already an entry on top of $a$, then the situation is precisely the situation described in the previous paragraph, shifted one row higher.  Therefore $F'(c) > F'(a)$ and $I(F'(c),F'(a))+I(F'(a),F'(b))-I(F'(c),F'(b))=1+0-0=1$ and hence $\{a,b,c\}$ is a type $B$ inversion triple.
\end{proof}

\begin{lemma}{\label{commute}}
The SSAF $\rho^{-1}(P)$ is the only SSAF with row entries $\{P_i \}$.
\end{lemma}

\begin{proof}
Let $K$ be an arbitrary SSAF whose row entries are given by the collection $\{P_i\}$.  Find the lowest row $r$ of $K$ whose entries appear in a different order from their appearance in row $r$ of $\rho^{-1}(P)=F'$.  Consider the largest element, $\delta$, of $P_r$ whose position in $K$ does not agree with its position in $F'$.  Then $\delta$ was placed in the cell on top of the left-most possible cell $a$ in $F'$ such that $\alpha=F'(a) \ge \delta$.  Therefore in $K$ the entry $\delta$ must reside in a cell $c$ of a column strictly to the right of the column containing $a$ in $K$.  The situation is depicted below, where $b$ is the cell immediately below the cell $c$ containing $\delta=K(c)$ in $K$ and the entry $\delta$ is shown as the entry in its cell.

\begin{center}
\begin{picture}(230,50)
\put(100,50){$K$}
\put(220,50){$F'$}
\put(0,0){ row $r-1$}
\put(0,20){ row $r$}
\put(80,15){$\tableau{ \\ a} \hdots \tableau{\delta \\ b}$}
\put(200,15){$\tableau{\delta \\ a} \hdots \tableau{ \\ b}$}
\end{picture}
\end{center}

If the column containing $\delta$ and the cell $b$ is strictly taller than the column containing the cell $a$ in $K$, then the triple $\{a,b,c\}$ is a type $B$ triple.  We have $I(\delta,K(a))+I(K(a),K(b))-I(\delta,K(b))=0+I(K(a),K(b))-0$.  We know that $K(a) < K(b)$ by Lemma \ref{typeb}.  So $I(K(a),K(b))$ must equal $0$ and hence $\{a,b,c\}$ is not a type $B$ inversion triple.

If the column containing $a$ is weakly taller in $K$ than the column containing $b$ and $c$, then the cell $d$ immediately on top of $a$ in $K$ contains an entry $K(d)$ which is less than $\delta$.  (If it were greater than $\delta$, this would contradict the assumption that $\delta$ is the greatest value for which the row placement differs.)  In particular, $K(d) < \delta=K(c)$ implies that $I(K(d),K(c))+I(K(c),K(a))-I(K(d),K(a))=0+0-0=0$.  So the triple $\{d,c,a\}$ is a type $A$ triple that is not a type $A$ inversion triple.

Therefore, regardless of which column is taller, $K$ contains at
least one non-inversion triple.  So $K$ is not a semi-skyline augmented filling.
\end{proof}

The above implies that the map $\rho^{-1}$ sends a reverse semi-standard Young tableau $P$ to the only possible semi-skyline augmented filling whose row entries are the same as the row entries of $P$.    In fact, let $\{R_i \}_j$ be any collection of sets of positive integers such that the number of entries in $R_i$ which are greater than or equal to $\alpha$ is great than or equal to the number of entries in $R_{i+1}$ which are greater than or equal to $\alpha$, for all $i \ge 1$ and all positive integers $\alpha$.  Any set with this property could be the set of row entries for a reverse semi-standard Young tableau.  So applying the map $\rho$ produces a semi-skyline augmented filling, the only semi-skyline augmented filling with row entries $\{R_i\}_j$.

\subsection{An analogue of Schensted insertion}{\label{insert}}

Schensted insertion is a procedure for inserting a positive integer $k$ into a semi-standard Young tableau $T$.  It is the fundamental operation of the Robinson-Schensted-Knuth (RSK) algorithm.  We define a procedure for inserting a positive integer $k$ into a semi-skyline augmented filling $F$ using a similar ``scanning and bumping" technique.  In Section {\ref{rsk_sec}} we use this procedure to describe an analogue of the RSK algorithm.

Let $F$ be a semi-skyline augmented filling of a weak composition $\gamma$ of $n$ into $m$ parts.  Then $F = (F(j))$, where $F(j)$ is the entry in the $j^{th}$ cell in reading order.  (Here we include the cells in the basement, so $j$ goes from $1$ to $n+m$.)  Let $\hat{j}$ be the cell immediately above $j$.  If the cell is empty, set $F(\hat{j})=0$.  We define the operation $k \rightarrow F$, for $k \le m$. 

\begin{procedure}{{\bf The insertion $k \rightarrow F$:}}

P1.  Set $i:=1$, set $x_1:=k$, set $p_0=\emptyset$, and set $j:=1$.

P2.  If $F(j) < x_i$ or $F(\hat{j}) \ge x_i$, then increase $j$ by $1$ and repeat this step.  Otherwise, set $x_{i+1}:=F(\hat{j})$ and set $F(\hat{j}):=x_i$.  Set $p_i=(a,b+1),$ where $(a,b)$ is the $j^{th}$ cell in reading order.  {\it (This means that the entry $x_i$ `bumps' the entry $x_{i+1}$ from the cell $p_i$.)}

P3.  If $x_{i+1} \not= 0$ then increase $i$ by $1$, increase $j$ by $1$, and repeat step P2.

P4.  Set $t_k$ equal to $p_i$ and terminate the algorithm.  

\end{procedure}

Notice that this procedure produces a sequence $I_k=(x_1, x_2, \hdots, x_r)$ of entries which are affected by the insertion of $k$, called the {\it insertion sequence}. The sequence $P_k=(p_1,p_2, \hdots, p_r)$ records which cells are affected by the insertion of $k$, with $t_k$ (the {\it termination cell}) denoting the new cell created at the termination of the insertion.  This sequence is called the {\it insertion path} of $k$.

\begin{lemma}{\label{higher_row}}
Let $x_i$ be an entry in $I_k$ and let $\alpha_i$ be the first occurrence in $F$ of the value $x_i$ after $p_{i-1}$ in reading order.  Then $x_i$ is inserted into a row above $\alpha_i$ by the procedure $k \rightarrow F$.
\end{lemma}

\begin{proof}
Let $s$ be the cell containing $\alpha_i$ in $F$ and let $a$ be the cell immediately above $s$.  The entry $x_i$ must appear before $a$ in the reading order of $F$ since $F$ is a non-attacking filling.  (If $x_i$ is $x_1$, we think of $x_i$ as appearing before all letters of $F$ in reading order.)  The entry in $a$ must be less than $x_i$ by the choice of $\alpha_i$.  Therefore if $x_i$ reaches $a$ during the insertion of $k$, then $x_i$ will replace the entry in $a$.  So $x_i$ appears in a row above $\alpha_i$.
\end{proof}

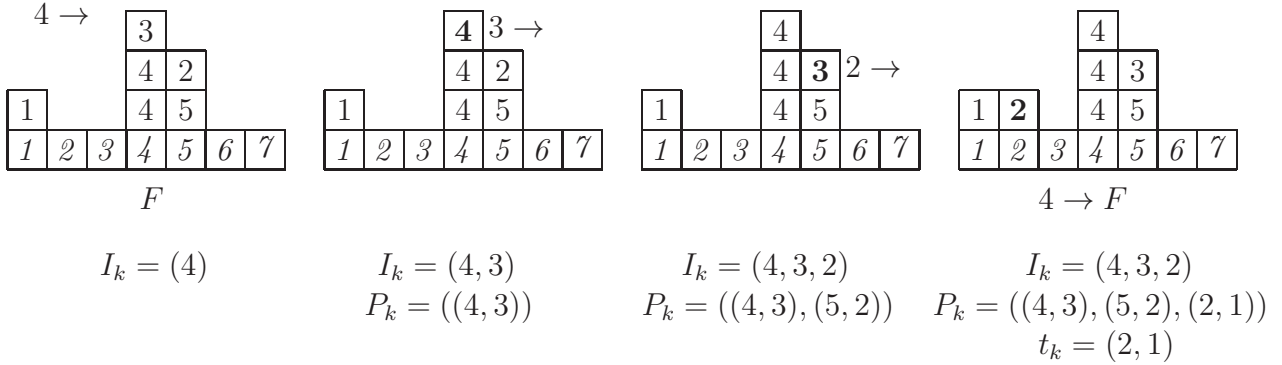
\begin{figure}
\begin{center}
\begin{picture}(465,120)
\put(50,50){$F$}
\put(0,110){$$\tableau{ & & & 3 \\ & & & 4 & 2 \\ 1 & & & 4 & 5 \\ {\it 1} & {\it 2} & {\it 3} & {\it 4} & {\it 5} & {\it 6} & {\it 7}}$$}
\put(10,120){$4 \rightarrow$}
\put(35,25){$I_k=(4)$}

\put(120,110){$$\tableau{ & & & {\bf 4} \\ & & & 4 & 2 \\ 1 & & & 4 & 5 \\ {\it 1} & {\it 2} & {\it 3} & {\it 4} & {\it 5} & {\it 6} & {\it 7}}$$}
\put(182,115){$3 \rightarrow$}
\put(140,25){$I_k=(4,3)$}
\put(135,10){$P_k=((4,3))$}

\put(240,110){$$\tableau{ & & & 4 \\ & & & 4 & {\bf 3} \\ 1 & & & 4 & 5 \\ {\it 1} & {\it 2} & {\it 3} & {\it 4} & {\it 5} & {\it 6} & {\it 7}}$$}
\put(317,100){$2 \rightarrow$}
\put(255,25){$I_k=(4,3,2)$}
\put(240,10){$P_k=((4,3),(5,2))$}

\put(360,110){$$\tableau{ & & & 4 \\ & & & 4 & 3 \\ 1 & {\bf 2} & & 4 & 5 \\ {\it 1} & {\it 2} & {\it 3} & {\it 4} & {\it 5} & {\it 6} & {\it 7}}$$}
\put(390,50){$4 \rightarrow F$}
\put(385,25){$I_k=(4,3,2)$}
\put(350,10){$P_k=((4,3),(5,2),(2,1) )$}
\put(390,-5){$t_k=(2,1)$}
\end{picture}
\end{center}
\caption{The insertion $4 \rightarrow F$}
{\label{insert4}}
\end{figure}

Lemma {\ref{higher_row}} implies that each entry $x_i$ is inserted into a row of $F$ which is higher than the next row containing an entry equal in value to $x_i$.  In particular, this means each entry appears in $k \rightarrow F$ above the basement entry of equal value, since every value appears in the basement.  Therefore the procedure $k \rightarrow F$ terminates in finitely many steps.  (See Figure \ref{insert4} for an example.)  To see that the resulting figure is indeed a semi-skyline augmented filling, we need the following lemma.

\begin{lemma}{\label{rowentries}}
The set of entries in a given row of the figure $k \rightarrow F$ is equal to the set of entries in the corresponding row of the figure $\insr$.
\end{lemma}

\begin{proof}
To prove this, we determine which entry of $F$ is bumped to a lower row of $F$ during $\insf$ and show that this is the same entry bumped to a lower row of $\rho(F)$ during $\insr$.

If $R_i$ is the collection of row entries appearing in row $i$ of $F$, let $\overline{R}_i = R_i \cup x_j$, where $x_j$ is the last entry bumped from a row higher than row $i$.  Let $|\overline{R}_i^{\ge u}|$ be the number of entries in $\overline{R}_i$ which are greater than or equal to $u$.  Similarly, let $|R_{i-1}^{\ge u}|$ be the number of entries in row $i-1$ of $F$ which are greater than or equal to $u$.  We claim that if $v$ is the largest entry in row $i$ such that $|\overline{R}_i^{\ge v}| > |R_{i-1}^{\ge v}|$, then $v$ is the entry in row $i$ which is bumped down to row $i-1$.  If no such $v$ exists, then no entry is bumped to row $i-1$ and the insertion procedure is terminated at row $i$.

To see this, let $v$ be the largest entry in $\overline{R}_i$ such that $|\overline{R}_i^{\ge v}| > |R_{i-1}^{\ge v}|$.  If no such $v$ exists, then the $m^{th}$ largest entry in $\overline{R}_i$ is less than or equal to $m$ entries in row $i-1$, for all $m$.  This means that each entry in $\overline{R}_i$ can be inserted into row $i$, and so the procedure terminates at row $i$, as predicted by the claim.

Assume next that such a $v$ exists.  Since $F$ is an SSAF, the pigeonhole principle implies that $|{R}_i^{\ge u}| \le |R_{i-1}^{\ge u}|$ for all $u$ in row $i$.  Therefore the largest $v$ which satisfies $|\overline{R}_i^{\ge v}| > |R_{i-1}^{\ge v}|$ must be less than or equal to $x_j$.  If $v=x_j$, then all the entries greater than or equal to $x_j$ in row $i-1$ are immediately below entries greater than $x_j$ and therefore $x_j$ is bumped down to the following row, as claimed.

If $v < x_j$, then there exists an entry in row $i-1$ of $F$ which is greater than or equal to $x_j$ but lies immediately beneath an entry, $\alpha$, less than $x_j$.  During the insertion, $x_j$ bumps the leftmost such $\alpha$.  Note that $\alpha$ must be greater than or equal to $v$ since all the entries greater than or equal to $v$ in row $i-1$ must lie beneath entries greater than or equal to $v$.  If $\alpha=v$, then all the entries greater than or equal to $\alpha$ in row $i-1$ lie immediately beneath entries greater than or equal to $\alpha$, so $\alpha$ is bumped to the next row as claimed.  If $\alpha>v$, there exists an entry smaller than $\alpha$ that sits on top of an entry greater than or equal to $\alpha$.  By the construction of $F$, this entry is to the right of $\alpha$, so $\alpha$ will bump this entry.  Again, since all entries greater than or equal to $v$ in row $i-1$ lie below entries greater than or equal to $v$, this entry is greater than or equal to $v$.  The same argument as above shows that either this bumped entry is $v$ or this entry bumps another entry in row $i$.  Therefore, the entry $v$ is eventually bumped.  Since all the entries greater than or equal to $v$ in row $i-1$ lie below entries greater than or equal to $v$, the entry $v$ cannot bump anything in row $i$.  Therefore $v$ is the entry bumped to the next row down, as claimed.

Now consider the insertion $\rho(F) \leftarrow k$.  This insertion procedure scans columns rather than rows, but the insertion path moves weakly south as it progresses through the columns from left to right.  To see this, consider an entry $\alpha$ bumped from row $i$ of $\rho(F)$.  Then $\alpha$ will be inserted into the column directly to its right on top of an entry weakly greater than $\alpha$.  The entry immediately to the right of $\alpha$ is strictly less than $\alpha$, so $\alpha$ cannot be inserted on top of this entry or at any higher row.  Therefore $\alpha$ is inserted into a row weakly lower than the row from which it is bumped.

If $\alpha$ is bumped to a row lower than row $i$, this means that the set of entries in row $i-1$ which are greater than or equal to $\alpha$ has cardinality less than the cardinality of the set containing $\alpha$ and the entries in row $i$ which are greater than or equal to $\alpha$.  Therefore an entry $\alpha$ is bumped to a lower row if $\alpha$ is the greatest entry in row $i$ such that $|\overline{R}_i^{\ge \alpha}| > |R_{i-1}^{\ge \alpha}|$, which is the same condition under which an entry is bumped to a lower row during $k \rightarrow F$.  Therefore the set of entries in a given row of the figure $k \rightarrow F$ is equal to the set of entries in the corresponding row of the figure $\insr$.
\end{proof}

\begin{proposition}{\label{produces_ssaf}}
If $F$ is a semi-skyline augmented filling and $k$ is an arbitrary positive integer, then the figure $k \rightarrow F$ is a semi-skyline augmented filling.  In particular, the insertion procedure commutes with the map $\rho$ in the sense that
$$\rho(k \rightarrow F ) = ( \rho(F) \leftarrow k),$$ where $\rho(F) \leftarrow k$ is the reverse Schensted insertion of $k$ into $\rho(F)$.
\end{proposition}

\begin{proof}
Lemma \ref{rowentries} states that the row entries of $k \rightarrow F$ are the same as those of $(\rho(F) \leftarrow k)$.  The semi-skyline augmented filling $\rho^{-1}(\rho(F) \leftarrow k)$ is the only semi-skyline augmented filling with these row entries by Lemma \ref{commute}, so we must prove that $\rho^{-1}(\rho(F) \leftarrow k)$ and $k \rightarrow F$ are the same semi-skyline augmented filling.  In other words, we must show that the row entries of $k \rightarrow F$ appear in the same positions as in $\rho^{-1}(\rho(F) \leftarrow k)$.

To see this, we prove that every triple in $k \rightarrow F$ is an inversion triple.  Consider first a type $A$ triple in $k \rightarrow F$, as depicted below.  

\begin{center}
\begin{picture}(50,40)
\put(0,15){$ \tableau{a \\ c} \hspace*{.5in} \tableau{b}$}
\end{picture}
\end{center}

Let $\alpha, \beta, \gamma$ be the entries contained in the cells $a,b,c$ respectively.  This triple is a non-inversion triple if $\alpha < \beta \le \gamma$.  The cell $a$ cannot be the termination cell for $k \rightarrow F$.  (If it were, then $c$ and $b$ would violate Lemma \ref{typeb} in $F$.)  If $b$ were the termination cell, then $\beta$ must have been bumped from a cell to the right of $a$ for otherwise $\beta$ would have bumped $\alpha$.  If the column containing $\beta$ before $\beta$ was bumped were taller than the column containing $a$, then $c$, the cell containing $\beta,$ and the cell below $\beta$ would form a type $B$ non-inversion triple in $F$.  Otherwise it must be that $a$, $c$, and the cell containing $\beta$ would form a type $A$ non-inversion triple in $F$.  Since all triples in $F$ are inversion triples, this is a contradiction.  Therefore a type $A$ non-inversion triple in $k \rightarrow F$ can only appear if some of the cells in the triple were bumped by the insertion.  

The entries in the cells of the triple are weakly increasing in reading order, so at most one of these entries is affected by the insertion procedure.  If $\alpha$ bumped a smaller entry, then the cells $a,b,$ and $c$ would form a type $A$ non-inversion triple in $F$.  The entry $\beta$ could not have bumped an entry from the cell $b$ for the same reason that it cannot be the termination cell.  Therefore the only cell in this triple that could be affected by the insertion procedure is $c$.  If $\gamma$ bumped a smaller entry, $\kappa$, from the cell $c$, then $\kappa < \beta$, for otherwise $a,b,c$ would be a type $A$ non-inversion triple in $F$.  

This implies that the entry $\kappa$ must be smaller than the entry $\delta$ in the cell $d$ directly below $b$, and so the entry $\epsilon$ in the cell $e$ directly below $c$ must be smaller than the entry $\delta$, for otherwise the cells $c,d,e$ form a type $A$ non-inversion triple in $F$.  But then if $\gamma$ appeared before $\beta$ in reading order, $\gamma$ would have bumped $\beta$.  So $\gamma$ was bumped from a cell after $b$ in reading order.  If $\gamma$ appeared on the same row of $F$ as $b$, the cell containing $\gamma$ would form an inversion triple with either $b$ and $d$ or $d$ and the cell immediately below $\gamma$ in $F$, depending on the column heights.  Therefore $\gamma$ must have been bumped from a cell on the same row as $c$.  Lemma \ref{typeb} implies that the column containing $\gamma$ in $F$ must be shorter than the column containing $a$, since $\gamma > \alpha$.  The entry above $\gamma$ must therefore be greater than $\beta$ so that these two cells together with the cell $b$ form a type $A$ inversion triple in $F$.  This means that $\gamma$ must be bumped by an entry greater than $\beta$.  Repeating the above argument implies that the entry which bumps $\gamma$ must come from the same row as $\gamma$, and continued repetition implies that every entry in the insertion must be bumped from the same row as $\gamma$, which contradicts the fact that the insertion algorithm begins at the highest cell in the diagram.  Therefore this situation cannot happen and every type $A$ triple is an inversion triple.

Next consider a type $B$ triple in $k \rightarrow F$, as depicted below.  

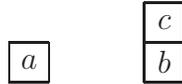
\begin{center}
\begin{picture}(40,40)
\put(0,15){$ \tableau{\\ a} \hspace*{.5in} \tableau{c \\ b}$}
\end{picture}
\end{center}

Let $\alpha, \beta, \gamma$ be the entries contained in the cells $a,b,c$ respectively.  This triple is a non-inversion triple if $\gamma \le \alpha < \beta$.  The cell $c$ cannot be the termination cell for $k \rightarrow F$.  To see this, notice that if $c$ were the termination cell, then $\gamma$ was bumped from a cell to the right of the column containing $a$.  (Otherwise $\gamma$ would be placed on top of $a$.)  But if $\gamma$ appears to the right of $a$ in $F$, then the cell containing $\gamma$, the cell beneath it, and the cell $a$ form a type $B$ non-inversion triple in $F$. 

If the insertion algorithm terminates at the cell $a$, then the entry $\alpha$ must have been bumped from a cell after $c$ in reading order.  If $\alpha$ were in the same row of $F$ as $c$, then Lemma \ref{typeb} implies that the column containing $\alpha$ is weakly shorter than the column containing $c$ and $b$, since $\alpha \le \beta$.  But then $c$, $b$, and the cell containing $\alpha$ form a type $A$ non-inversion triple.  If $\alpha$ were in the same row of $F$ as $b$, then $c,b,$ and the cell containing $\alpha$ in $F$ would constitute a type $B$ non-inverstion triple unless the column containing $\alpha$ were weakly taller in $F$ than the column containing $b$ and $c$.  In this case, the entry which bumps $\alpha$ must be less than $\beta$ since it must be inserted on top of an entry that is greater than $\beta$.  (The entry below $\alpha$ in $F$ is greater than $\beta$ since these two cells together with $b$ form a type $A$ inversion triple in $F$.)  Then we are in the same situation as before and can repeat the argument indefinitely to show that no entry is bumped from the row above $a$, which contradicts the fact that the insertion algorithm begins at the first cell in reading order.  Therefore a type $B$ non-inversion triple in $k \rightarrow F$ can only appear if some of the cells in the triple were bumped by the insertion.  

The entries in the cells of the triple are weakly decreasing in reading order, so at most one of these entries is affected by the insertion procedure.  If $\gamma$ bumped a smaller entry from the cell $c$, then this entry would be less than $\alpha$ and Lemma \ref{typeb} implies that the triple $a,b,c$ would be a type $B$ non-inversion triple in $F$.  If $\alpha$ bumps an entry from the cell $a$, then $\alpha$ must be bumped from a cell after $c$ in reading order, for otherwise $\alpha$ would bump $\gamma$.  If the entry $\alpha$ is bumped from the row containing $c$, then the column containing $\alpha$ in $F$ must be weakly shorter than the column containing $c$ and $b$ since $\alpha > \beta$.  But then the cell containing $\alpha$ together with $c$ and $b$ form a type $A$ non-inversion triple.  Therefore $\alpha$ must be bumped from the row containing $a$.  Lemma \ref{typeb} implies that the column containing $\alpha$ in $F$ must be taller than or equal in height to the column containing $c$, since $\gamma \le \alpha$.  The entry above $\alpha$ must therefore be greater than $\gamma$ so that these two cells together with the cell $c$ form a type $A$ inversion triple in $F$.  This means that $\alpha$ must be bumped by an entry greater than or equal to $\gamma$.  Repeating the above argument implies that the entry which bumps $\alpha$ must come from the same row as $\alpha$, and continued repetition implies that every entry in the insertion must be bumped from the same row as $\alpha$, which contradicts the fact that the insertion algorithm begins at the highest cell in the diagram.  Therefore the entry $\alpha$ could not have bumped a smaller entry from the cell $a$.  The entry $\beta$ could not have bumped a smaller entry from the cell $b$, for if so then $a,b,c$ would be a type $B$ non-inversion triple in $F$by Lemma \ref{typeb}.  Therefore this situation cannot happen and every type $B$ triple is an inversion triple.
\end{proof}

The insertion procedure $k \rightarrow F$ is closely connected to Schensted insertion.  Let {\it reverse Schensted insertion} be the variation of Schensted insertion which maps a positive integer $k$ into a reverse semi-standard Young tableau~\cite{Stanley2}.  This is equivalent to altering the Schensted insertion procedure described by Knuth \cite{Knu70} by reversing the directions of the inequalities, and scanning columns rather than rows.

\subsection{The bijection $\Psi$ between SSYT and SSAF}{\label{bijection}}

Let $T$ be a semi-standard Young tableau.  We may associate to $T$ the word $\mathit{col}(T)$, which consists of the entries from each column of $T$, read top to bottom from columns left to right, as in Figure {\ref{columns}}.  In general, any word $w$ can be decomposed into its maximal strictly decreasing subwords, called {\it column words}.  This decomposition is called $\mathit{col}(w)$.   For example, if $w= 3 \; 5 \; 4 \; 2 \; 2 \; 1$, then $\mathit{col}(w) = 3 \cdot 5 \; 4 \; 2 \cdot 2 \; 1$.

\begin{figure}
$$\tableau{10 \\ 9 & 11 \\ 8 & 10 & 10 \\ 4 & 7 & 8 \\ 2 & 5 & 5 \\ 1 & 2 & 3 & 5 & 10}$$
\caption{Here $\mathit{col}(T) = $ $10$ $9$ $8$ $4$ $2$ $1$ $\cdot$ $11$ $10$ $7$ $5$ $2$ $\cdot$ $10$ $8$ $5$ $3$ $\cdot$ $5$ $\cdot$ $10$}
{\label{columns}}
\end{figure}

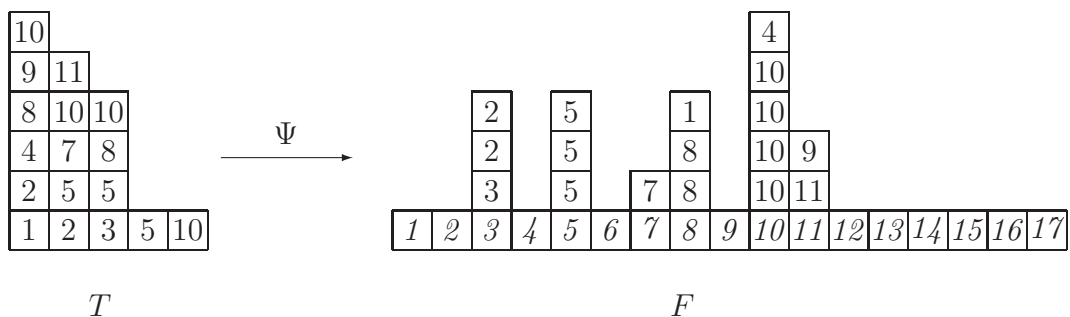
\begin{figure}[b]
\begin{center}
\begin{picture}(400,120)
\put(0,100){$$\tableau{10 \\ 9 & 11 \\ 8 & 10 & 10 \\ 4 & 7 & 8 \\ 2 & 5 & 5 \\ 1 & 2 & 3 & 5 & 10}$$}
\put(30,0){$T$}

\put(80,60){\vector(1,0){50}}
\put(100,65){$\Psi$}

\put(145,100){$$\tableau{ & & & & & & & & & 4 \\ & & & & & & & & & 10 \\ & & 2 & & 5 & & & 1 & & 10 \\ & & 2 & & 5 & & & 8 & & 10 & 9 \\ & & 3 & & 5 & & 7 & 8 & & 10 & 11 \\ {\it 1} & {\it 2} & {\it 3} & {\it 4} & {\it 5} & {\it 6} & {\it 7} & {\it 8} & {\it 9} & {\it 10} & {\it 11} & {\it 12} & {\it 13} & {\it 14} & {\it 15} & {\it 16} & {\it 17}}$$}
\put(250,0){$F$}
\end{picture}
\end{center}
\caption{The map $\Psi:T \rightarrow F$ takes a SSYT $T$ to an SSAF $F$.}
\label{mapping}
\end{figure}

Begin with an arbitrary SSYT $T$ and the empty SSAF $\phi$ whose basement row contains all the letters of $\mathbb{Z}_+$.  Let $k$ be the rightmost letter in $\mathit{col}(T)$.  Insert $k$ into $\phi$ to obtain the SSAF $(k \rightarrow \phi)$.  Then let $k'$ be the next letter in $\mathit{col}(T)$ reading right to left.  Obtain the SSAF $(k' \rightarrow (k \rightarrow \phi) )$.  Continue in this manner until all the letters of $\mathit{col}(T)$ have been inserted.  The resulting diagram is the SSAF $\Psi(T)$.  (See Figure {\ref{mapping}}.)

Lemma {\ref{rowentries}} implies that that the entries in the rows of $\Psi(T)$ are precisely the entries in the rows of the reverse SSYT $P$ obtained by applying reverse Schensted insertion to $col(T)$.  The shape of $P$ is equal to the transpose of the shape of $T$ (\cite{Knu70}, \cite{RSK}).  This implies that the shape of $\Psi(T)$ is a rearrangement of the shape of $T$.

The map $\Psi$ is therefore a weight-preserving map from semi-standard Young tableaux of shape $\lambda$ to semi-skyline augmented fillings whose shape rearranges $\lambda$.  We now show that $\Psi$ is a bijection.  

Note that $\Psi$ is invertible, since reverse Schensted insertion is an invertible procedure and $\Psi$ commutes with reverse Schensted insertion by Proposition {\ref{produces_ssaf}}.  To describe the inverse directly, consider an arbitrary SSAF $F$.  Consider the set $S$ of nonzero columns of $F$.  (In Figure {\ref{mapping}}, $S$ contains the third, fifth, seventh, eighth, tenth, and eleventh columns.)  The highest cell in each of these columns is a cell which was created during the reverse Schensted insertion of the leftmost column of $col(T)$, since the SSAF produced by the insertion of all columns except the leftmost column has the shape of the diagram of $T$ with the leftmost column omitted.  Pick the shortest such column, $C_1$, where if two columns have equal height then the column farther to the right is considered to be shorter.

The cell $c_1$ at the top of column $C_1$ of $F$ is the cell at which the insertion of the leftmost letter of $col(T)$ terminates.  Delete $F(c_1)$ and scan the reading word in reverse order to determine which entry (if any) bumped the entry $F(c_1)$.  (Any entry, $\beta$, greater than $F(c_1)$ beneath an entry less than $F(c_1)$ would have bumped $F(c_1)$.)  If such an entry is found, replace it with $F(c_1)$ and then scan the reading word in reverse order beginning at that position to determine the previous bumping entry and replace it with the entry which bumped $F(c_1)$.  Continue until the first letter of the reading word is reached.  The last bumping entry is therefore the first entry of $col(T)$.  

For example, in Figure {\ref{mapping}}, the entry $7$ in the seventh column was the last to be placed into $F$.  It was bumped by the $9$ in the eleventh column, which was bumped by the $10$ in the fourth row of the tenth column.  This $10$ was not bumped from any other position, so the $10$ is the first letter in the word $col(T)$.  Removing the $7$ in this manner produces the SSAF pictured in Figure {\ref{newssaf}}.

\begin{figure}
$$\tableau{ & & & & & & & & & 4 \\ & & & & & & & & & 9 \\ & & 2 & & 5 & & & 1 & & 10 \\ & & 2 & & 5 & & & 8 & & 10 & 7 \\ & & 3 & & 5 & & & 8 & & 10 & 11 \\ {\it 1} & {\it 2} & {\it 3} & {\it 4} & {\it 5} & {\it 6} & {\it 7} & {\it 8} & {\it 9} & {\it 10} & {\it 11} & {\it 12} & {\it 13} & {\it 14} & {\it 15} & {\it 16} & {\it 17}}$$
\put(250,0){$F$}
\caption{The SSAF obtained by removing the $7$ in the seventh column of $F$.}
\label{newssaf}
\end{figure}

Repeat this procedure using the resulting SSAF for the second smallest column $C_2$ in $S$ to obtain the second letter of $col(T)$.  Continue in this manner until the highest cell in each of the columns in $S$ has been removed.  At this point the first column of $T$ has been recovered.  Then repeat the entire procedure with the new SSAF, whose nonzero columns are each one cell shorter than the nonzero columns in the original SSAF.  Continue this procedure until there are no columns remaining.  The resulting word is $col(T)$.

The map $\Psi:SSYT \rightarrow SSAF$ is a weight-preserving, shape-rearranging bijection.  This proves combinatorially that $$s_{\lambda}(x)=\sum_{\gamma^+=\lambda} \widehat{E}_{\gamma}(x;0,0),$$ since the left-hand side consists of all monomials $x^T$ where $T$ is a semi-standard Young tableau of shape $\lambda$ while the right-hand side consists of all monomials $x^F$ where $F$ is a semi-skyline augmented filling of shape $\gamma$, where $\gamma^+=\lambda$.

\subsection{A basis for $\mathbb{Q}[x]$}

Let $\gamma$ be a weak composition of $n$ into $m$ parts, where $m \in \mathbb{Z}^+ \cup \{ \infty \}$.  The monomial corresponding to the composition $\gamma$ is $x^{\gamma}=\prod_{i=1}^m x_i^{\gamma_i}$.  The monomials corresponding to compositions with infinitely many parts form a $\mathbb{Q}$-basis for $\mathbb{Q}[x]$, where $x=\{x_1, x_2, \hdots \}$.

Recall that the reverse dominance order on weak compositions $\gamma$ and $\mu$ is given by

$$\mu \le \gamma \iff \sum_{i=k}^{\infty}\mu_{i} \le \sum_{i=k}^{\infty}\gamma_{i}, \hspace*{.3in} \hspace*{.1in}  \forall k \ge 1.$$

Let $NK_{\gamma,\mu}$ be the coefficient of $x^{\mu}$ in $\widehat{E}_{\gamma}(x;0,0)$.  Then $NK_{\gamma,\mu}$ is equal to the number of semi-skyline augmented fillings of shape $\gamma$ and content $\mu$.  The bijection presented in Section~{\ref{bijection}} implies that $$K_{\lambda,\mu} = \sum_{\gamma^+=\lambda} NK_{\gamma,\mu},$$ where $K_{\lambda, \mu}$ is the Kostka number, equal to the number of SSYT of shape $\lambda$ and weight $\mu$.

\begin{proposition} {\label{basis2}} Suppose that $\gamma$ and $\mu$ are both weak compositions of $n$ into $k$ parts and $NK_{\gamma, \mu} \not= 0$.  Then $\mu \le \gamma$ in the reverse dominance order.  Moreover, $NK_{\gamma,\gamma} = 1$.
\end{proposition}

\begin{proof}
Assume that $NK_{\gamma, \mu} \not=0$.  By definition, this means that there exists a semi-skyline augmented filling of shape $\gamma$ and content $\mu$.  Suppose that an entry $m$ appears in one of the first $m-1$ columns.  Then this column contains a descent, since the basement entry of this column is less than $m$.  Therefore all entries greater than or equal to $m$ must appear after the $(m-1)^{th}$ column.  So $\sum_{i=m}^{\infty} \mu_i \le \sum_{i=m}^{\infty} \gamma_i$ for each $m$, as desired.  

Moreover if $\mu=\gamma$, then the entries in the $i^{th}$ column must all be equal to $i$.  To see that this augmented filling is a semi-skyline augmented filling, consider first a type $A$ triple as shown.
\begin{center}
\begin{picture}(60,30)
\put(0,15){$$\tableau{\alpha \\ \alpha}$$}
\put(45,15){$$\tableau{\delta}$$}
\put(29,20){.}
\put(22,20){.}
\put(36,20){.}
\end{picture}
\end{center}
Here $\delta > \alpha$, so $I(\alpha, \delta)+I(\delta,\alpha)-I(\alpha,\alpha)=0+1-0=1$ implies that the triple is indeed a type $A$ inversion triple.

Next consider a type $B$ triple as shown below.
\begin{center}
\begin{picture}(60,30)
\put(0,0){$$\tableau{\alpha}$$}
\put(45,15){$$\tableau{\beta \\ \beta}$$}
\put(21,4){.}
\put(36,4){.}
\put(28,4){.}
\end{picture}
\end{center}
Since $\alpha < \beta$, we have $I(\beta, \alpha)+I(\alpha,\beta)-I(\beta,\beta)=1+0-0=1$.  Therefore the triple is indeed a type $B$ inversion triple.

Since all triples are inversion triples and there are no descents, the filling is a semi-skyline augmented filling.  It is the only SSAF of shape $\gamma$ and content $\gamma$, so $NK_{\gamma,\gamma}=1$.
\end{proof}

\begin{cor}
The polynomials $\widehat{E}_{\gamma}(x;0,0)$ are a $\mathbb{Q}$-basis for $\mathbb{Q}[x]$.
\end{cor}

\begin{proof}
Proposition {\ref{basis2}} is equivalent to the assertion that the transition matrix from the polynomials $\widehat{E}_{\gamma}(x; 0,0)$ to the monomials (with respect to the reverse dominance order) is upper triangular with 1's on the main diagonal.  Since this matrix is invertible, the $\widehat{E}_{\gamma}(x;0,0)$ are a basis for all polynomials.
\end{proof}

Notice that restriction to compositions of $n$ into $m$ parts implies that the $\widehat{E}_{\gamma}(x; 0,0)$ form a basis for polynomials of degree $n$ in $m$ variables.

\section{An analogue of the Robinson-Schensted-Knuth Algorithm}{\label{rsk_sec}}

The insertion procedure utilized in the above bijection is an analogue of Schensted insertion, the fundamental operation of the Robinson-Schensted-Knuth (RSK) Algorithm.

\begin{theorem}{(Robinson-Schensted-Knuth~\cite{RSK})}
There exists a bijection between $\mathbb{N}$-matrices of finite support and pairs of semi-standard Young tableaux of the same shape.
\end{theorem}

We apply the same procedure to arrive at an analogue of the RSK Algorithm for semi-skyline augmented fillings.  Recall that Theorem {\ref{rsk}} states that there exists a bijection between $\mathbb{N}$-matrices of finite support and pairs of semi-skyline augmented fillings whose shapes are rearrangements of the same partition.  Composing the RSK algorithm with the bijection $\Phi$ between semi-standard Young tableaux and semi-skyline augmented fillings proves Theorem {\ref{rsk}}.  We now provide a direct bijection that commutes with this composition.

\subsection{The map $\Phi: \mathbb{N}$-matrices $\longrightarrow$ SSAF $\times$ SSAF}

Let $A=(a_{i,j})$ be an $\mathbb{N}$-matrix with finite support.  There exists a unique two-line array corresponding to $A$ which is defined by the non-zero entries in $A$.  Let $a_{i,j}$ be the first non-zero entry encountered when scanning the entries of $A$ from left to right, top to bottom.  Place an $i$ in the top line and a $j$ in the bottom line $a_{i,j}$ times.  When this has been done for each non-zero entry, one obtains the following array. 
\[ w_A= \left( \begin{array}{ccc}
i_1 & i_2 & . . . \\
j_1 & j_2 & . . .
\end{array} \right) \]

Notice that if $i_r=i_{r+1}$, then $j_r \le j_{r+1}$ in $w_A$.  Every two-line array with this property can be obtained from a matrix with finite support.  We may therefore consider two-line arrays with this property instead of matrices of finite support.

\begin{procedure}{{\bf The map $\Phi: w_A \longrightarrow \text{SSAF} \times \text{SSAF}$:}}

P1.  Set $r:=l$, where $l$ is the length of $w_A$.  Let $F=\phi=G$, where $\phi$ is the empty SSAF.  

P2.  Set $F:=(j_r \rightarrow F)$.  Let $h_r$ be the height of the column in $(j_r \rightarrow F)$ at which the insertion procedure $(j_r \rightarrow F)$ terminates.

P3.  Place $i_r$ on top of the leftmost column of height $h_r-1$ in $G$ such that doing so does not create a descent.  Set $G$ equal to the resulting figure.

P4.  If $r-1 \not= 0$, repeat step P2 for $r:=r-1$.  Else terminate the algorithm.

\end{procedure}

Notice that the entries in the top row of the array are weakly increasing when read from left to right.  This means that if $h_r > 1$, placing $i_r$ on top of the leftmost column of height $h_r-1$ in $G$ does not create a descent.  If $h_r =1$, we claim that the $i_r^{th}$ column of $G$ does not contain an entry from a previous step.  To see this, argue by contradiction.  Assume there exists an entry $i_m$ in the $i_r^{th}$ column, where $m > r$.   Then the entries $j_m, j_{m-1}, \hdots j_r$ of the array $w_A$ form a weakly decreasing sequence by the nature of the array.  The properties of the reverse Schensted insertion of $j_r$ imply that the termination of the insertion of $j_r$ into $F$ must occur in a higher row of $F$ than the termination of the insertion of $j_m$.  This contradicts the assumption that $i_r$ and $i_m$ are inserted into columns of the same height.

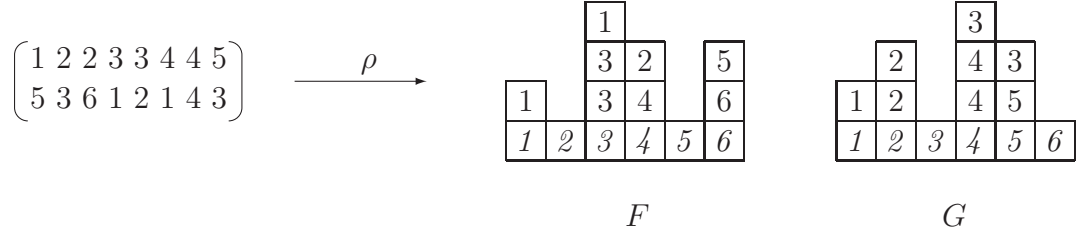
\begin{figure}
\begin{picture}(360,90)
\put(4,55){\oval(10,30)[l]}
\put(5,45){5  3  6  1  2  1  4  3}
\put(5,60){1  2  2  3  3  4  4  5}
\put(80,55){\oval(10,30)[r]}
\put(105,55){\vector(1,0){50}}
\put(130,60){$\rho$}
\put(185,70){$$\tableau{ & & 1 \\ & & 3 & 2 & & 5 \\ 1 & & 3 & 4 & & 6 \\ {\it 1} & {\it 2} & {\it 3} & {\it 4} & {\it 5} & {\it 6}}$$}
\put(230,0){$F$}
\put(310,70){$$\tableau{ & & & 3 \\ & 2 & & 4 & 3 \\ 1 & 2 & & 4 & 5 \\ {\it 1} & {\it 2} & {\it 3} & {\it 4} & {\it 5} & {\it 6}}$$}
\put(350,0){$G$}
\end{picture}
\caption{$\Phi: A \rightarrow SSAF \times SSAF$}
{\label{fg}}
\end{figure}

In this way the shape of $G$ becomes a rearrangement of the shape of $F$.  When the process is complete, the result is a pair $(F,G)=\Phi(A)$ of fillings whose shapes are rearrangements of the same partition (see Figure {\ref{fg}}).  To see that $F$ and $G$ are indeed semi-skyline augmented fillings, we prove that the map $\Phi$ produces precisely the pair of SSAFs obtained by applying the reverse RSK algorithm to $A$ followed by the map $\rho^{-1}$.

\begin{proposition}
{\label{result}}  If $(P,Q)$ is the pair of reverse semi-standard Young tableaux obtained by applying the reverse RSK algorithm to the matrix $A$, then $(\rho^{-1}(P),\rho^{-1}(Q))=\Phi(A)=(F,G)$.
\end{proposition}

\begin{proof}
Proposition {\ref{produces_ssaf}} implies that $\rho^{-1}(P)=F$.  We must prove that $\rho^{-1}(Q)=G$, where $Q$ is the recording tableau.  Consider the insertion of an element $i_r$ into $G$.  The position of $i_r$ in $G$ is determined by the termination of insertion $(j_r \rightarrow F)$.  Let $h_r$ denote the height of the column at which the insertion of $(j_r \rightarrow F)$ terminates.  The entry $i_r$ appears in row $h_r$ of $G$ by definition.  By Proposition {\ref{produces_ssaf}}, the reverse Schensted insertion of $j_r$ terminates in row $h_r$ of $P$.  Hence $i_r$ appears in row $h_r$ of $Q$.  Therefore the entries in the rows of $Q$ are the same as those in the rows of $G$.

Recall that $\rho^{-1}(Q)$ is the unique SSAF whose rows contain precisely the row entries of $Q$.  The row entries of $Q$ are inserted into $\rho^{-1}(Q)$  in decreasing order, so that an entry $\alpha$ appears in the leftmost available column such that the entry immediately below $\alpha$ in this column is greater than or equal to $\alpha$.

The top line of the array $w_A$ is a weakly increasing sequence, so its entries are inserted into $G$ in weakly decreasing order.  Therefore the entries in a given row $r$ are inserted into $G$ in decreasing order.  An entry $i_r$ is placed on top of the leftmost possible column which does not create a descent.  To see that the entries in a given row of $\rho^{-1}(Q)$ appear in the same columns as the entries in a given column of $G$, argue by induction on the entries in the top line of the array $w_A$.

Consider first the entry $i_l$, where $l$ is the length of $w_A$.  This entry appears in the first row of $G$, in the $i_l^{th}$ column.  This entry also appears in the $i_l^{th}$ column of $\rho^{-1}(Q)$, since that is the leftmost possible position for $i_l$ in $\rho^{-1}(Q)$. 

Assume that the entries $i_l, i_{l-1}, \hdots, i_r$ appear in the same columns in $G$ as they do in $\rho^{-1}(Q)$.  Consider the entry $i_{r-1}$.  This entry appears on top of the entry $\alpha$ in $\rho^{-1}(Q)$.  Therefore $\alpha \ge i_{r-1}$.  This fact and the fact that $\alpha$ appears in a lower row of $G$ than $i_{r-1}$ imply that $\alpha$ is inserted into $G$ before $i_{r-1}$.  Since all the entries which were inserted into $G$ before $i_{r-1}$ appear in the same columns as in $\rho^{-1}(Q)$, there is no entry on top of $\alpha$ in $G$ during the insertion of $i_{r-1}$.  This means that if $i_{r-1}$ is not placed on top of $\alpha$, then $i_{r-1}$ is placed in a column to the left of the column containing $\alpha$.  Let $\beta$ be the entry under $i_{r-1}$ in this situation.  Then $\beta$ appears in this same position in $\rho^{-1}(Q)$.  Since $i_{r-1}$ is not on top of $\beta$ in $\rho^{-1}(Q)$, there must be an entry $\gamma > i_{r-1}$ on top of $\beta$ in $\rho^{-1}(Q)$.  But if this were the case, then $\gamma$ would also be on top of $\beta$ in $G$ since by assumption, all entries greater than $i_{r-1}$ are in the same positions in $G$ as they are in $\rho^{-1}(Q)$.  This contradicts the assumption that $i_{r-1}$ is placed on top of $\beta$.  Therefore the entries of $G$ are in precisely the same positions as in $\rho^{-1}(Q)$, so $\rho^{-1}(Q)=G$.
\end{proof}

Proposition {\ref{result}} implies in particular that $F$ and $G$ are indeed semi-skyline augmented fillings.  The semi-skyline augmented filling $G$ records the column to which a cell is added in $F$.  Therefore $F$ and $G$ are rearrangements of the same shape.  These two facts imply that $\rho$ is a map from $\mathbb{N}$-matrices to pairs of semi-skyline augmented fillings which are rearrangements of the same shape. 

\subsection{The inverse of the map $\Phi$}

Let $(F,G)$ be a pair of semi-skyline augmented fillings whose shapes are rearrangements of the same partition $\lambda$.  Then $(\rho(F),\rho(G))$ is a pair of reverse semi-standard Young tableaux of shape $\lambda$.  Since the reverse Robinson-Schensted-Knuth algorithm is invertible, this pair $(\rho(F),\rho(G))$ comes from a unique matrix $A$ of finite support.  This matrix is $\Phi^{-1}((F,G))$.

We describe the computation of this inverse without mapping through the reverse semi-standard Young tableau pair.  Let $G_{rs}$ be the highest occurrence of the smallest entry of $G$ in reading order.  (Here $G_{rs}$ is the element of $G$ in row $r$ and column $s$.)  

Let $s'$ be the rightmost column of height $r$ in $F$.  Then $F_{rs'}$ is the termination cell for the insertion of the last letter, $j_1$ in the array corresponding to the matrix $A$.  Delete $F_{rs'}$ from $F$ and $G_{rs}$ from $G$.  Scan right to left, bottom to top (backwards through the reading word) starting with the cell directly to the left of $F_{rs'}$ to determine which cell (if any) bumped $F_{rs'}$.  If there exists a cell $k$ before $F_{rs'}$ in the reading word such that $F(k) > F_{rs'}$ and the cell directly on top of $k$ has value less than or equal to $F_{rs'}$, this $F(k)$ bumped $F_{rs'}$ by the argument from Section {\ref{bijection}}.  Replace $F(k)$ by $F_{rs'}$ and repeat the procedure with $F(k)$ starting from the cell $k$.  Continue working backward through the reading word until there are no more letters in the reading word.  The resulting entry is the letter $j_1$.  (Notice that this is the same procedure used in the map $\Psi^{-1}$ in Section {\ref{bijection}}, so we have already proved that this procedure does in fact yield the inverse of the insertion of $j_1$.)

Next find the highest occurrence of the smallest entry $j_2$ of the new recording filling.  Repeat the procedure to find $i_2$.  Continue until there are no more entries in the insertion and recording tableaux.  Then all of the $i$ and $j$ values of the array $w_A$ have been determined.

\begin{cor}
The RSK algorithm commutes with the above analogue $\Phi$.  That is, if $(P,Q)$ is the pair of SSYT produced by the RSK algorithm applied to a matrix $A$, then $(\Psi(P),\Psi(Q))=\Phi(A)$.
\end{cor}

\begin{proof}
The pair $(P,Q)$ maps bijectively to the pair $(S,T)$ of reverse SSYT produced by applying the reverse RSK algorithm to $A$.  Proposition {\ref{result}} states that $(\rho^{-1}(S),\rho^{-1}(T))=(F,G)=\Phi(A).$  Proposition {\ref{produces_ssaf}} implies that $(\Psi(P),\Psi(Q))=(\rho^{-1}(S),\rho^{-1}(T))$, so $(\Psi(P),\Psi(Q))=(F,G)$ as desired.
\end{proof}

The algorithm $\Phi$ retains many of the properties of the Robinson-Schensted-Knuth Algorithm.  Several of these properties are listed below, although this list is by no means complete.

\begin{cor}{\label{symm}}
Assume that $\Phi(A)=(F,G)$.  Then $\Phi(A^t)=(G,F)$.
\end{cor}

\begin{cor}{\label{knuth}}
Let $\Phi(A)=(F,G)$.  Let $B$ be a matrix such that the bottom line of $w_B$ is Knuth equivalent to the bottom line of $w_A$.  Then $\Phi(B)=(F,G')$ for some SSAF $G'$.
\end{cor}

Corollaries {\ref{symm}} and {\ref{knuth}} follow directly from the fact that the algorithm $\Phi$ commutes with the Robinson-Schensted-Knuth Algorithm.

\subsection{Standardization}{\label{standard}}

A standard Young tableau of shape $\lambda$ is a semi-standard Young tableau whose weight is $\prod_{i=1}^n x_i$, where $n=|\lambda|$.   Any semi-standard Young tableau $T$ can be mapped to a standard Young tableau through a procedure called standardization.

\begin{definition}
Let $\gamma$ be a weak composition of $n$.  A {\it skyline augmented filling (SAF)} of shape $\gamma$ is an SSAF of weight $\prod_{i=1}^n x_i$.  Let $F_{\gamma}$ denote the number of skyline augmented fillings of shape $\gamma$.
\end{definition}

Notice that $F_{\gamma}=0$ for certain compositions $\gamma$.  This means that our analogue of standardization, {\it skylining}, must alter the shape of such compositions.  

\begin{procedure}[$Skylining(F)$]
Assume that $F$ is an SSAF of shape $\gamma$, where $| \gamma |=n$.  Consider the reading word as a collection of row words for $F$.  Standardize the reading word for $F$ in the usual manner.  Then place the new entries back into the original shape of $F$, maintaining the same order.  If the lowest entry, $\alpha_j$ of the $j^{th}$ column is not equal to $j$, shift the entire column to rest immediately on top of the $\alpha_j^{th}$ column.  Denote the result by $sk(F)$.
\end{procedure}

The entries in the first row of $F$ are always strictly increasing from left to right. Therefore the shifting does not permute the order of the non-zero columns.   Notice that if $a$ appears before $c$ in reading order, then $F(a) \le F(c) \Leftrightarrow sk(F(a)) < sk(F(c))$.  Therefore $I(F(a),F(b))+I(F(b),F(c))-I(F(a),F(c))$ remains the same after skylining.  The same is true for $I(F(c),F(a))+I(F(a),F(b))-I(F(c),F(b))$ since $F(a) < F(c)$.  This means that skylining preserves the inversion triples.  This fact and the fact that there are no descents within the columns of $F$ imply that the figure $sk(F)$ is indeed a semi-skyline augmented filling.

\begin{figure}
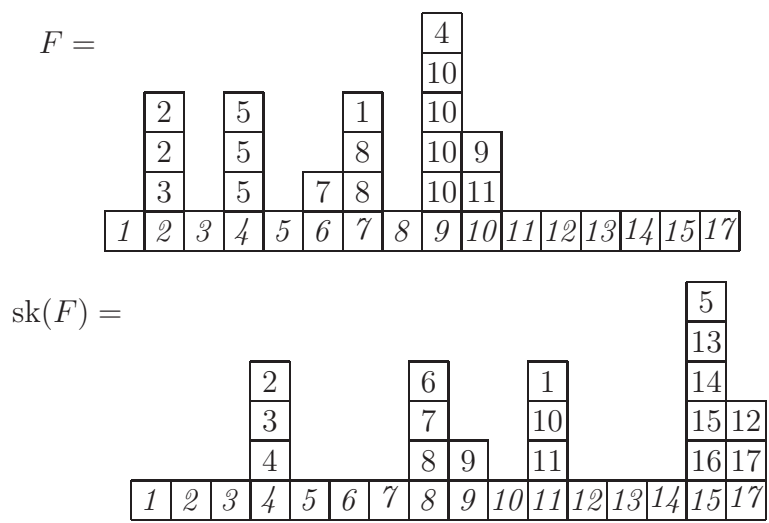
{\label{skylining}}
$$F=\tableau{ & & & & & & & & 4\\ & & & & & & & & 10\\ & 2 & & 5 & & & 1 & & 10 \\
 & 2 & & 5 & & & 8 & & 10 & 9 \\  & 3 & & 5 & & 7 & 8 & & 10 & 11 \\ 
{\it 1} & {\it 2} & {\it 3} & {\it 4} & {\it 5} & {\it 6} & {\it 7} & {\it 8} & {\it 9} & {\it 10} & {\it 11} & {\it 12} & {\it 13} & {\it 14} & {\it 15} & {\it 17}} $$

$$\text{sk}(F)=\tableau{ & & & & & & & & & & & & & & 5\\ & & & & & & & & & & & & & & 13\\ & & & 2 & & & & 6 & & & 1 & & & & 14 \\
& & & 3 & & & & 7 & & & 10 & & & & 15 & 12 \\  & & & 4 & & & & 8 & 9 &  &11 & & & & 16 & 17 \\ 
{\it 1} & {\it 2} & {\it 3} & {\it 4} & {\it 5} & {\it 6} & {\it 7} & {\it 8} & {\it 9} & {\it 10} & {\it 11} & {\it 12} & {\it 13} & {\it 14} & {\it 15} & {\it 17}} $$
 \caption{An SSAF $F$ and its standardization $\text{sk}(F)$}
 \end{figure}

\begin{proposition}
The bijection $\Psi$ commutes with standardization.
\end{proposition}

\begin{proof}
Let $T$ be a semi-standard Young tableau and $\Psi(T)$ the corresponding SSAF.  We saw in section {\ref{rev}} that the SSAFs are in bijective correspondence with reverse SSYTs, so it is enough to consider the reverse SSYT $R$ corresponding to $T$.

The set of entries in a given row $r$ of $R$ are sent to a new set of row entries under standardization.  However, the new entries are the same entries which appear in row $r$ of the skylining of $\Psi(T)$ by the definition of $sk(\Psi(T))$.  There is only one SSAF with the same row entries as the standardization of $R$, so this must be $sk(\Psi(T))$.
\end{proof}

\bibliographystyle{plain}

\end{document}